\documentclass[12pt,fleqn]{amsart}

\usepackage{amssymb}
\usepackage{enumerate}

\usepackage{tikz-cd}

\newtheorem{theorem}{Theorem}[section]
\newtheorem{lemma}[theorem]{Lemma}
\newtheorem{corollary}[theorem]{Corollary}
\newtheorem{proposition}[theorem]{Proposition}

\newtheorem{problem}[theorem]{Problem}
\newtheorem{example}[theorem]{Example}

\theoremstyle{definition}
\newtheorem{construction}[theorem]{Construction}
\newtheorem{definition}[theorem]{Definition}
\newtheorem{remark}[theorem]{Remark}

\newcommand{\propR}{\mathcal{(R)}}
\newcommand{\propE}{\mathcal{(E)}}

\newenvironment{scclaim}[1][Claim.]
    {\medskip \noindent \textsc{#1} }
    {\medskip}

\renewcommand{\phi}{\varphi}
\newcommand*\diff{\mathop{}\!\mathrm{d}}

\newcommand{\eps}{\varepsilon}
\newcommand{\fA}{\mathfrak A}

\newcommand{\fB}{\mathfrak B}

\newcommand{\cA}{{\mathcal A}}

\newcommand{\cM}{{\mathcal M}}
\newcommand{\cN}{{\mathcal N}}
\newcommand{\cI}{{\mathcal I}}

\newcommand{\cC}{{\mathcal C}}

\newcommand{\cU}{{\mathcal U}}
\newcommand{\cV}{{\mathcal V}}

\newcommand{\btu}{\bigtriangleup}

\newcommand{\ult}{{\rm ult}}
\newcommand{\vf}{\varphi}

\newcommand{\lra}{\longrightarrow}
\newcommand{\sm}{\setminus}

\newcommand{\sub}{\subseteq}

\newcommand{\cde}{\protect{\rm CDE} }
\newcommand{\fin}{\protect{\rm fin} }

\newcommand{\wt}{\widetilde}
\newcommand{\idE}{\boldsymbol{\cI}}
\newcommand{\nonE}{\protect{\rm non}(\idE)}

\title{Countable discrete extensions of compact lines}
\author[M. Korpalski]{Maciej Korpalski}
\author[G.\ Plebanek]{Grzegorz Plebanek}
\address{Instytut Matematyczny, Uniwersytet Wroc\l awski, pl.\ Grunwaldzki 2, 50-384 Wroc\-\l aw\\ Poland}
\email{Maciej.Korpalski@math.uni.wroc.pl \\ Grzegorz.Plebanek@math.uni.wroc.pl}
\date{}

\thanks{The first author was partially supported by the grant 2018/29/B/ST1/00223 from National Science Centre, Poland.
The second author was supported by the NCN (National Science Centre, Poland), under the Weave-UNISONO call in the Weave programme 2021/03/Y/ST1/00124.
}

\begin{document}

\begin{abstract}
We consider a separable compact line $K$ and its extension $L$ consisting of $K$
and a countable number of isolated points.
The main object of study is the existence of a bounded extension operator
$E: C(K)\to C(L)$. We show that if such an operator exists then there is
one for which $\|E\|$ is an odd natural number.
We prove that if the topological weight of $K$ is bigger than or equal to the least
cardinality of a set $X\sub [0,1]$ than cannot be covered by a sequence of closed sets
of measure zero then there is an extension $L$ of $K$ admitting no bounded extension operator.
\end{abstract}
\maketitle

\section{Introduction}

A \textit{compact line} is a compact space whose topology is defined by means of a linear order.
The double arrow space
\[ \mathbb{S} = \big((0, 1] \times \{0\}\big) \cup \big([0,1)\times \{1\}\big),\]
ordered lexicographically is a familiar example of a compact line.
The space $\mathbb{S}$ is nonmetrizable, but separable and first countable, see \cite[Exercise 3.10.C]{En89}.
There is a natural generalization of the double arrow space:
we can consider an arbitrary closed subset $F$ of the unit interval, any set $X \subseteq F$ and define
\[F_X = \big(F \times \{0\}\big) \cup \big( X \times \{1\}\big).\]
As before, the space $F_X$, ordered lexicographically, is a separable compact line and it is nonmetrizable whenever $X$ is uncountable.
In fact, for an infinite set $X$, the space $F_X$ is of topological weight $|X|$; see \cite{Ma08} for more information.
It turns out that spaces of the form $F_X$ exhaust the class of separable compact lines:

\begin{theorem}[Ostaszewski, \cite{Os74}] \label{int:1}
The space $K$ is a separable compact linearly ordered space if and only if $K$ is homeomorphic to $F_X$ for some closed set $F \subseteq [0, 1]$ and a subset $X\subseteq F$.
\end{theorem}

Given a pair of compact spaces $K\sub L$, by an \textit{extension operator} $E:C(K)\to C(L)$ we mean a bounded linear operator (between Banach spaces of continuous functions) such that $Ef|K = f$ for every $f\in C(K)$.
By the classical Borsuk-Dugundji extension theorem such an operator of norm one exist whenever $K$ is metrizable, see \cite{Pe68} or \cite[II.4.14]{LT73}.
Following \cite{CL65} we write $\eta(K,L)$ for the infimum of the norms of all possible extension operators $E:C(K)\to C(L)$ if there are any; thus $\eta(K,L)=\infty$ means that there is no bounded extension operator.
There are many examples of nonmetrizable compacta $K\sub L$ with $\eta(K,L)<\infty$; calculating $\eta(K,L)$ may be, however, quite involved, see e.g. \cite{AM15}.

If $K$ is any compact space, we call a superspace $L \supseteq K$ a \textit{countable discrete extension} of $K$ and write $L\in\cde(K)$ if $L$ is compact and the set $L\sm K$ is a countable infinite discrete set.
The main subject of this paper is investigating two properties of countable discrete extensions of separable compact lines.
In a sense, those properties measure the complexity of the way in which isolated points are added to the initial space.

\begin{definition} \label{def:properties}
Given a fixed compact space $K$ and $L\in\cde(K)$,
\begin{enumerate}[(i)]
\item $L$ has property $\propR$ if there is a continuous retraction from $L$ onto $K$;
\item $L$ has property $\propE$ if $\eta(K,L)<\infty$.
\end{enumerate}
\end{definition}

Properties $\propR$ and $\propE$ were considered in a series of papers
\cite{DP17}, \cite{MP18}, \cite{AMP20} in connection with twisted sums of Banach spaces, see the final section for more information.
We shall try to convince the reader that the subject of our study is quite subtle and interesting on its own.
Recall also that there is a vivid trend of investigating properties of  Banach spaces of continuous functions on compact lines, see e.g. \cite{Ma08}, \cite{KK12}, \cite{CT14} and \cite{Mi17}.

We show here the following results concerning a separable compact line $K$ of topological weight $w(K)$
(in fact,  the separability is not needed for first two results).

\begin{enumerate}[(a)]
  \item For every $L\in\cde(K)$, either $L$ has property $\propR$ (so $\eta(K,L)=1$) or $\eta(K,L)\ge 3$.
  \item For every $L\in\cde(K)$, if $\eta(K,L)<5$ then $\eta(K,L)\le 3$. We also sketch tyhe argument for proving
  that $\eta(K,L)$ is either infinite or equal to an odd natural number.
  \item If $w(K)\ge\omega_1$ then there is $L\in\cde(K)$ such that $\eta(K,L)=3$.
  \item If $w(K)\ge\nonE$ then there is $L\in\cde(K)$ without property $\propE$.
\end{enumerate}

Here $\nonE$ denotes the least cardinality of a set $X\sub [0,1]$ that cannot be covered by a countable family of closed null sets.
For constructing examples mentioned above we find it convenient to see zero-dimensional compact lines as Stone spaces of algebras generated by almost chains in some countable set.
We do not know if it is relatively consistent that (for a separable $K$ as above), if $w(K)=\omega_1$ then every countable discrete extension $L$ of $K$ has property $\propE$.
As it is explained in section \ref{final}, such a fact is true for every {\em nonseparable} compact line.

We wish to thank Witold Marciszewski for very fruitful conversations concerning the subject; in particular, we present in section \ref{between} his result communicated to us privately.

\section{Preliminaries} \label{sec:Preliminaries}

All topological spaces we consider are Hausdorff. In the sequel, $K$ and $L$ always stand for compact spaces.
For a space $K$ we denote by $C(K)$ the Banach space of continuous functions $f: K \to \mathbb{R}$ equipped with the supremum norm denoted by $\|\cdot\|$.
By the Riesz representation theorem we can identify the dual space $C(K)^\ast$ with the space of signed regular Borel measures of bounded variation which is denoted by $M(K)$.
For a measure $\mu \in M(K)$ we can use the Jordan decomposition theorem to write $\mu = \mu^+ - \mu^-$ for some non-negative orthogonal measures $\mu^+, \mu^-$.
The total variation of $\mu$ is then $\|\mu\| = \mu^+(K) + \mu^-(K)$. Given $x\in K$, $\delta_x$ is the probability Dirac measure concentrated at the point $x$.

In the sequel, we write $\omega$ for the set of natural numbers equipped with the discrete topology.
We often consider countable discrete extensions of a given compact space $K$ in the form $K\cup\omega$, tacitly assuming that $K\cap\omega=\emptyset$.
As declared above, we write $L \in \cde(K)$ to denote that $L$ is a countable discrete extension of $K$.

\begin{remark}
If $r: L \to K$ is a continuous retraction, then $C(K) \ni f \to f \circ r \in C(L)$ defines an extension operator of norm one.
This proves that every countable discrete extension with property $\propR$ has property $\propE$.
\end{remark}

Extensions having property $\propR$ are, in a sense, trivial.
Most of easy constructions of countable discrete extensions have this property.
Recall that there are spaces which do not have extensions without property $\propR$.

\begin{example} \label{ex:Metrizable has R}
If a space $K$ is metrizable, then every $L \in \cde(K)$ has property $\propR$.
\end{example}

\begin{proof}
Since $K$ is compact and metrizable, it is separable
\cite[Theorem 4.1.18]{En89} and has a countable base. It follows that the space $L$ also
has a countable base and therefore it is metrizable \cite[Theorem 4.2.8]{En89}.
Let us fix a compatible metric $d$ on $L$.
We can define a retraction $r$ by mapping each element in $L$ to the closest element in $K$ (to any of them if there are many). Such a function is well defined as $K$ is compact, so the closest element always exists.

Checking that $r$ is continuous amounts to verifying that whenever a sequence
$(x_n)_{n \in \omega}$ in $\omega$ converges to $x \in K$, then $r(x_n) \to r(x) = x$.

Given $\eps > 0$ we have $d(x_n,x)<\eps$ for almost all $n$ so
\[d(r(x_n), x_n) \leq d(x, x_n) < \eps \text{ and } d(r(x_n), x) < 2\eps,\]
for large $n$, by the triangle inequality.
\end{proof}

We recall below useful characterizations of properties $\propR$ and $\propE$, see \cite[Lemma 2.7]{MP18} for the proof which is fairly standard; cf.\ \cite{Pe68} and \cite{Ve71}.

\begin{lemma} \label{p:1}
Let $K$ be any compact space and let $L \in \cde(K)$.
\begin{enumerate}[(a)]
  \item $L$ has property $\propR$ if and only if there is a sequence of points $(x_n)_{n \in \omega}$ in $K$ such that for every function $f \in C(L)$ we have
  \[ \lim_{n \to \infty} (f(x_n) - f(n)) = 0.\]
  \item $L$ has property $\propE$ if and only if there is a bounded sequence of signed measures $(\mu _n)_{n \in \omega}$ on $K$ such that
   $\mu_n-\delta_n\to 0 $ in the $weak^\ast$ topology of $C(L)^\ast$, i.e.
    for every $f \in C(L)$ we have
  \[ \lim_{n \to \infty} \Big(\int_K f \diff \mu_n - f(n)\Big) = 0.\]
\end{enumerate}
\end{lemma}

\begin{remark} \label{p:2}
Concerning Lemma \ref{p:1}(b), the norm of the extension operator $E$ satisfies
$\|E\| = \sup_{n\in\omega} \|\mu_n\|$.
\end{remark}

Recall that there are spaces of arbitrarily large weight that do not admit countable discrete extensions without property $\propR$.
Indeed, take any cardinal number $\kappa$ and consider the Cantor cube $2^\kappa$.
Then $2^\kappa$ is an absolute retract in the class of compact zero-dimensional spaces so, in particular, every $L\in CDE(2^\kappa)$ has property $\propR$. This can be demonstrated directly as follows.

The space $2^\kappa$ has a subbase consisting of sets
\[C_\alpha^i = \{x\in 2^\kappa : x(\alpha) = i\},\]
for $\alpha < \kappa$ and $i=0,1$.
For every $\alpha$, $L$ can be partitioned into clopen sets $\widetilde{C_\alpha^i}$ such that
$\widetilde{C_\alpha^i}\cap K=C_\alpha^i$.
Thus we can define a continuous retraction $r$ by $r|K = id_K$ and, for $n \in L \sm 2^\kappa$,
let $r(n)$ be the only point in
\[ \bigcap_{\alpha \in \kappa, i\in \{0,1\}} \{C_\alpha^i : n \in \widetilde{ C_\alpha^i}\}.\]

\section{Calculating $\eta(K,L)$}\label{calculating}

Corson and Lindenstrauss \cite{CL65} showed that if $K$ is a one point compactification of an uncountable
discrete space then for every compact superspace $L\supseteq K$, if $\eta(K.L)<\infty$
then $\eta(K,L)$ is an odd natural number.
We show in this section that the same phenomenon is present in our context.

Throughout this section we assume that $K$ is an arbitrary, not necessarily separable, compact line and $L=K\cup\omega$ is its countable discrete extension. We denote simply by $<$ the linear order on $K$ and
write, for instance, $[s,t]=\{x\in K: s\le x\le t\}$ for $s,t\in K$. Denote by $\alpha$ and $\beta$ the least and the greatest element in $K$, respectively.

We first give a very technical but convenient criterion for convergence of measures on $L=K\cup\omega$.

\begin{lemma}\label{ca:0}
Let $(\nu_n)_n$ be a bounded sequence in $M(K)$ such that $\nu_n(K)=1$ for every $n$.
Suppose that whenever $s,t\in K$, $s<t$, and closed subsets $F,H$ of $L=K\cup\omega$ satisfy 
$F\cap K\sub [\alpha,s]$ and $H\cap K\sub [t,\beta]$ then

\begin{enumerate}[(i)]
\item $\nu_n[t,\beta]=0$ for almost all $n\in F\cap\omega$;
\item $\nu_n[\alpha,s]=0$ for almost all $n\in H\cap\omega$.
\end{enumerate}

Then $\nu_n-\delta_n\to 0$ in the $weak^\ast$ topology of $M(L)$.
\end{lemma}

\begin{proof}
Recall first that if $(t_n)_{n\in\omega}$ is a sequence in a compact topological space $T$ and
$\cU$ is any non-principal ultrafilter on $\omega$ then there is a unique element $T\ni t=\lim_{n\to\cU} t_n$ such that $\{n\in\omega: t_n\in V\}\in\cU$ for every open set $V$ containing $t$.

To prove the lemma suppose that the assertion does not hold; then the sequence of measures $\nu_n-\delta_n$ has a nonzero cluster point $\mu$.
Take an ultrafilter $\cU$ such that
\[\mu=\lim_{n\to\cU} (\nu_n-\delta_n)=\lim_{n\to\cU}\nu_n- \lim_{n\to\cU}\delta_n\neq 0. \]
Of course, $\lim_{n\to\cU}\delta_n=\delta_s$ for some $s\in K$, so writing $\nu=\lim_{n\to\cU}\nu_n$ we have $\nu\neq \delta_s$.
Since $\nu(K)=1$, we conclude that either $|\nu|[\alpha,s)>0$ or $|\nu|(s,\beta]>0$;
suppose, for instance, that the latter holds (the former case will follow by a symmetric argument).

Then there is $K\ni t_1>s$ such that $|\nu|[t_1,\beta]>0$.
Apply the normality of $L$ to the sets $[\alpha,s] \sub L\sm [t_1,\beta]$ to find an open set $V\sub L$ such that
\[ [\alpha,s]\sub V\sub\overline{V}\sub L\sm [t_1,\beta],\]
and put $t_0=\sup (\overline{V}\cap K)$. Note that $t_0 < t_1$.
Now, by $(i)$, for every $t \in [t_1,\beta]$ we have $\nu_n[t,\beta] \to 0$ for $n\in V$ (note that the set $V\cap\omega$ is infinite).

Observe  that any $g\in C(K)$ vanishing on $[\alpha, t_1]$ can be uniformly approximated by step functions
built on intervals $[t,t')$ contained in $(t_1,\beta]$ so we conclude that $\int_K g \diff\nu_n\to 0$
for $n\in V$ which yields $|\nu|[t_1,\beta]=0$, a contradiction.
\end{proof}

\begin{theorem}\label{ca:1}
 If $\eta(K,L)<3$, then $L$ has property $\propR$.
 \end{theorem}

\begin{proof}
By Lemma \ref{p:1} and Remark \ref{p:2} there is a sequence of measures $\mu_n$ on $K$ such that $c=\sup_n \|\mu_n\|<3$ and $\mu_n-\delta_n\to 0$ in the $weak^\ast$ topology.
Fix $\delta>0$ such that $c+3\delta<3$ and for any $x\in K$ set
\[ A_x=\{n\in\omega: \mu_n^+[\alpha,x]\ge 1-\delta\}.\]
Then for every $n\in\omega$ we define
\[ x_n=\inf\{x\in K: n\in A_x\}.\]
Note that $\mu_n(K)\to 1$, so $x_n$ is well-defined for almost all $n\in\omega$.

\begin{scclaim}
Consider $s,t\in K$ with $s<t$; let $F,H$ be closed subsets of $L$ such that $F\cap K \sub [\alpha, s]$ and $H\cap K\sub [t,\beta]$. Then the following sets are finite
\[ I=\{n\in F: x_n\ge t\}, \quad J=\{n\in H: x_n\le s\}.\]
\end{scclaim}

To check Claim consider a continuous function $f:L\to [0,1]$ such that $f(x)=1$ for $x\le s$ and $f(x)=0$ for $x\ge t$, and a function $g:L\to [0,1]$ defined as $g=1-f$.

If we suppose that $I$ is infinite then for $n\in I$ we have $\mu_n^+[\alpha,x]<1-\delta$ whenever $x<t$ so
\[ \int_K f \diff \mu_n\le  \int_K f \diff \mu_n^+ < 1-\delta;\]
on the other hand, $\lim_{n\in I} f(n)=1$, a contradiction with $\mu_n-\delta_n\to 0$.

Suppose now that $J$ is infinite. Then, as $\lim_{n \in J} \int_K g \diff \mu_n - g(n) = 0$ and $g$ is equal to $1$ on $H$, we have
\[\int_K g \diff \mu_n \ge 1-\delta\mbox{ and } \int_K f \diff \mu_n < \delta,\]
for almost all $n\in J$.
At the same time, $\mu_n^+[\alpha, s]\ge 1-\delta$, so examining the integral on the right hand side above, $\mu_n^-[\alpha, s]\ge 1-\delta$ must eventually hold for $n\in J$.
It follows that
\[ |\mu_n|(K)\ge |\mu_n|[\alpha,s]+|\mu_n|[t,\beta]\ge 2(1-\delta)+1-\delta=3-3\delta>c,\]
contrary to our assumption $\|\mu_n\|\le c$.
\medskip

Once we have verified Claim, we conclude from Lemma \ref{ca:0} that $f(x_n)-f(n)\to 0$ for every $f\in C(L)$, and we are done.
\end{proof}

\begin{theorem}\label{ca:2}
If $\eta(K,L)<5$, then there is an extension operator $E:C(K)\to C(L)$ with $\|E\|\le 3$.
\end{theorem}

\begin{proof}
By Lemma \ref{p:1} and Remark \ref{p:2} there is a sequence of measures $\mu_n$ on $K$ such that $c=\sup_n \|\mu_n\|<5$ and $\mu_n-\delta_n\to 0$ in the $weak^\ast$ topology.
Fix $\delta>0$ such that $c+3\delta<5$.

We shall define a sequence of measures $\nu_n$ of norm at most $3$ satisfying 
$\mu_n-\delta_n\xrightarrow{weak^\ast} 0$. For the rest of the proof we assume that $\|\mu_n\|>3-\delta/4$ for every $n$; for other measures we can just put $\nu_n = \mu_n$.
For any $x\in K$ we define
\begin{align*}
A_x^0 &=\{n\in\omega: \mu_n^+[\alpha,x] \ge 1-\delta/4\}, &&x_n^0 =\inf\{ x\in K: n\in A_x^0\}; & & \\
A_x^1 &=\{n\in A_x^0: \mu_n^-[\alpha,x] \ge 1-\delta/2\}, &&x_n^1 =\inf\{ x\in K: n\in A_x^1\}; & & \\
A_x^2 &=\{n\in A_x^1: \mu_n^+[\alpha,x] \ge 2-\delta\}, &&x_n^2 =\inf\{ x\in K: n\in A_x^2\}. & &
\end{align*}
Note that, since $|\mu_n|(K)>3-\delta/4$ for every $n$ and $\mu_n(K)\to 1$, we can and do assume that all three sets on the right hand side above are nonempty for all $n$ so $x_n^i$ are well-defined and $x_n^0 \le x_n^1 \le x_n^2$.

We consider the sequence of measures
\[ \nu_n=\delta_{x_n^0}- \delta_{x_n^1}+\delta_{x_n^2},\]
and prove that $\nu_n-\delta_n \xrightarrow{weak^\ast} 0$ referring to Lemma \ref{ca:0}.
Fix $s,t\in K$ with $s<t$, and suppose that $F$ and $H$ are closed sets in $L$ such that $F\cap K \sub [\alpha, s]$ and $H\cap K\sub [t,\beta]$.
Take a continuous function $f:L\to [0,1]$ such that $f|[\alpha,s]=1$ and $f|[t,\beta]=0$; also set the function $g=1-f$.
Now let us check the assumptions of Lemma \ref{ca:0} in a few steps.

\begin{scclaim}[Step 1:]
The set $I=\{n\in F \cap\omega: s<t\le x_n^0\le x_n^1\le x_n^2\}$ is finite.
\end{scclaim}

We know that $\mu_n-\delta_n\xrightarrow{w^\ast} 0$ and $\lim_{n\in F} f(n) = 1$ (if $F \cap \omega$ is infinite), so the statement follows from the fact that $\lim_{n\in I} \int_K f\diff\mu_n=1$ for infinite $I\sub F$.

\begin{scclaim}[Step 2:]
The set $I=\{n\in H\cap\omega: x_n^0\le  s < t\le x_n^1\le x_n^2\}$ is finite.
\end{scclaim}

If $I$ were infinite then we would have $\lim_{n\in I} \int_K f\diff\mu_n=0$, while for $n\in I$
\[ \int_K f\diff\mu_n\ge \int_{[\alpha,s]} f\diff\mu_n^+ - \int_{[\alpha,t]} f\diff\mu_n^-\ge
1-\delta/4+\delta/2-1=\delta/2.\]

\begin{scclaim}[Step 3:]
The set $I=\{n\in F\cap\omega : x_n^0\le x_n^1\le s < t\le x_n^2\}$ is finite.
\end{scclaim}

Indeed, for infinite $I$ we would have $\lim_{n\in I} \int_K f\diff\mu_n=1$, while for $n\in I$
\[ \int_K f \diff\mu_n\le \int_{[\alpha,t]} f\diff\mu_n^+ - \int_{[\alpha,s]} f\diff\mu_n^-\le
2-\delta+\delta/2-1=1-\delta/2.\]

\begin{scclaim}[Step 4:]
The set $I=\{n\in H\cap\omega: x_n^0\le x_n^1\le x_n^2\le s < t\}$ is finite.
\end{scclaim}

For infinite $I$ we would again have $\lim_{n\in I} \int_K f\diff\mu_n=0$.
This, together with $x_n^2\le s$ meaning $\mu_n^+[\alpha,s]\ge 2-\delta$, implies that $\mu_n^-[\alpha,s]\ge 2-\delta$ for large $n\in I$.
Consequently, $|\mu_n|[\alpha,s]\ge 4-2\delta$ must hold eventually for $n\in I$.
On the other hand, $\lim_{n\in I} \int_K g\diff\mu_n=1$ implies $|\mu_n|[s,\beta]\ge 1-\delta$ for almost all $n\in I$ and we get a contradiction with $\|\mu_n\|\le c<5-3\delta$.

\begin{scclaim}[Step 5:]
Other cases that would violate the assumption of Lemma \ref{ca:0} are also excluded.
\end{scclaim}
\vspace{-12pt}

If we suppose, for instance, that the set $I=\{n\in F \cap\omega:  x_n^0\le s\le x_n^1\le t\le x_n^2\}$ is infinite then we can split it into two parts and apply on of the above cases.
\medskip

Now, by Lemma \ref{ca:0}, $\nu_n-\delta_n\xrightarrow{weak^\ast} 0$, as required.
\end{proof}

Examining the proofs of Theorem \ref{ca:1} and Theorem \ref{ca:2}, one can conclude that the argument may be further generalized.
Checking the details seems to be tedious, so we only sketch the general idea.

\begin{theorem}\label{ca:3}
If $\eta(K,L)<\infty$, then $\eta(K,L)$ is an odd natural number.
\end{theorem}

\begin{proof}
For the proof one can consider a natural number $k$ such that $2k-1\le \eta(K,L)<2k+1$ and $\delta>0$ that
is small enough. Proceeding by induction we can assume that $\|\mu_n\|\ge 2k-1-\delta$, where
the measures $\mu_n$ are related to an extension operator of norm $<2k+1$.
Then we can define for $x\in K$ the sets  $A_x^0,\ldots, A_x^{2k-2}$ by
\begin{align*}
A_x^0 &=\{n\in\omega : \mu_n^+[\alpha,x]\ge 1-\delta/2^{k}\}, \\
A_x^{2j+1} &=\{n\in A_x^{2j} : \mu_n^-[\alpha,x]\ge j+1-\delta/2^{k-2(j+1)}\},\\
A_x^{2j} &=\{n\in A_x^{2j-1} : \mu_n^+[\alpha,x]\ge j+1-\delta/2^{k-2j}\}.
\end{align*}
Then, after setting $x_n^i=\inf\{x\in K: n\in A_x^i\}$ for $i=0,\dots, 2k-2$, we consider the measures
\[ \nu_n=\sum_{i=0}^{2k-2} (-1)^i\delta_{x_n^i};\]
clearly, $\|\nu_n\|\le 2k-1$ so it remains to check that $\nu_n-\delta_n\xrightarrow{weak^\ast} 0$.
\end{proof}

\section{Countable discrete extensions and Boolean algebras} \label{cde}

In this section we  describe a method of constructing countable discrete extensions of separable compact lines via Stone spaces of Boolean algebras of subsets of $\omega$.
We  use here the classical Stone duality, referring to  \cite{Ko89} if necessary.
Given a Boolean algebra $\fA$,  $\ult(\fA)$ denotes  the Stone space of ultrafilters on $\fA$ which is a compact zero-dimensional space with a base consisting of all clopens of the form
\[\widehat{A} = \{\cU\in \ult(\fA) : A\in \cU\}, \]
for $A\in\fA$.

The basic idea is simple:
If an algebra $\fA\subseteq P(\omega)$ contains $\fin$, the ideal of finite subsets of $\omega$, then $\ult(\fA)$ becomes a compactification of $\omega$, by identifying principal ultrafilters with natural numbers.
Hence, once we can represent our basic compact zero-dimensional space $K$ as $K=\ult(\fA/\fin)$ for some Boolean algebra $\fA$ of subsets of $\omega$ (or any other countable set) then $L=\ult(\fA)$ is a countable discrete extension of $K=\ult(\fA/\fin)$.
Our first objective is to understand properties $\propR$ and $\propE$ in the Boolean language.

\begin{lemma} \label{cde:1}
Suppose that $K = \ult(\fA/\fin)$ for some algebra $\fA$ of subsets of $\omega$ containing $\fin$.
Then $L = \ult(\fA)$ has property $\propR$ if and only if there exists a lifting $\theta: \fA/\fin \to \fA$.
\end{lemma}

Here $\theta : \fA/\fin \to \fA$ is said to be a \textit{lifting} if it is a Boolean algebra
monomorphism such that $\pi~ \circ~ \theta =~id_{\fA/\fin}$.
The proof of the above lemma is standard, see e.\ g.\ \cite{DP17}.

To give an analogous lemma on $\propE$ recall that if $K=\ult(\fB)$ for any Boolean algebra $\fB\sub P(\omega)$ then we can think that $M(K)$ is identified with $M(\fB)$, the space of signed finitely additive measures on $\fB$ having bounded variation.
In this setting the norm of a measure $\mu \in M(\fA)$ is given by
$\|\mu\|=|\mu|(\omega)$, where the variation $|\mu|$ is defined for $A\in\fB$ as
\[|\mu|(A)=\sup_{B\in \fA, \, B\sub A} \big(|\mu(B)| + |\mu(A\setminus B)|\big).\]

\begin{lemma} \label{cde:2}
In the setting of Lemma \ref{cde:1},
 $L\in \cde(K)$ has property $\propE$ if and only if there is
 a bounded sequence $(\mu_n)_{n \in \omega}$ in $M(\fA)$
 such that
 \begin{enumerate}[(i)]
 \item $\mu_n(I)=0$ for every $I\in\fin$ and every $n$;
 \item
$ \lim_{n \to \infty} (\mu_n(A) - \delta_n(A)) = 0$ for every $A\in\fA$.
\end{enumerate}
\end{lemma}

\begin{proof}
This follows from Lemma \ref{p:1}(b) and the following observations.

There is an obvious correspondence between finitely additive measures on
$\fA/\fin$ and finitely additive measures on $\fA$ vanishing on finite sets.
Note also that for any zero-dimensional compact space $L$ and a sequence $\nu_n$ in $M(L)$, $\nu_n\to 0$ in the $weak^\ast$ topology if and only if $\nu_n(C)\to 0$ for every clopen $C\sub L$;
we can apply this remark to $\nu_n=\mu_n-\delta_n\in M(\fA)$.
\end{proof}

Recall that for $A,B\sub\omega$, $A\sub^\ast B$ stands for the relation of almost inclusion and means that  the set $A\sm B$ is finite;
likewise $A=^\ast B$ means that the set $A\btu B $ is finite.

It is well-known that, in terms of Stone duality, zero-dimensional  compact lines correspond to chain algebras; a chain algebra is one having linearly ordered set of generators.
The following Lemma is essentially known, see \cite[Theorem 8.7]{MP18}.

\begin{lemma}\label{cde:3}
Let $X \sub [0, 1]$ and suppose that $\cA=\{A_x:x\in X\}$ is an almost chain of subsets of $\omega$, that is $A_x\sub^\ast A_y$ whenever $x,y\in X$ and $x<y$.
Denote by $\fA$ the Boolean algebra generated by $\cA\cup\fin$.

Then $K=\ult(\fA/\fin)$ is a separable compact line with $w(K) = |X|$ and $L=\ult(\fA)$ is a countable discrete extension of $K$.
\end{lemma}

\begin{proof}
We have already explained that  $L$ may be seen as  a countable discrete extension of $K = \ult(\fA/\fin)$.
Then $K$ is a compact line, as $\fA/\fin$ is generated by a chain,
see e.g.\ \cite[Theorem 15.7]{Ko89}.
Recall that this follows from the fact that every ultrafilter $\cU\in\ult(\fA/\fin)$ is uniquely determined by the set $X(\cU)=\{x\in X: A_x/\fin\in \cU\}$ and we can suitably order $\ult(\fA/\fin)$ by declaring $\cU\le \cV$ when $X(\cV)\sub X(\cU)$.

Finally, $K$ is separable: take a countable set $D\sub X$ such that
for every $x\in X$ and $\delta>0$ there is $d\in D$ such that $x-\delta<d\le x$.
If for any $x\in X$ we denote by $\cU_x$ the unique ultrafilter in $K$
such that $x$ is the first element in $X(\cU)$, then $\{\cU_d: d\in D\}$ is easily seen to be dense in $K$.
\end{proof}

We can also reverse this characterization in the following manner.

\begin{lemma}\label{cde:4}
Let $K$ be a zero-dimensional separable compact line and let $L\in\cde(K)$.
Then there is $X\sub [0, 1]$ and an almost chain $\cA=\{A_x:x\in X\}$ of subsets of some countable set $N$ such that
\begin{enumerate}[(i)]
    \item $K$ is homeomorphic to $\ult(\fA/\fin)$ and
    \item $L$ is homeomorphic to $\ult(\fA)$,
\end{enumerate}
where $\fA$ is the algebra generated by $\cA\cup\fin(N)$.
\end{lemma}

\begin{proof}
By Theorem \ref{int:1} we know that $K$ is homeomorphic to  the space $F_X$ for some closed set $F \sub [0,1]$ and a subset $X \sub F$, so for the proof we consider $K=F_X$.
Note that, as $K$ is zero-dimensional, $X$ must be dense in $F$ (with respect to the natural topology).
As $L\in\cde(F_X)$, we have $L=F_X\cup N$ for some countable infinite set $N$ of isolated points.

For every $x\in X$, the set
\[  C_x=\{(y,i) \in F \times \{0,1\}: (y, i) <_{lex} (x, 1)\}  \]
is clopen in $F_X$ so there is a clopen set $\widetilde{C_x}$ in $L$ such that $\widetilde{C_x}\cap F_X=C_x$.
Consider $A_x=\widetilde{C_x}\cap N$.

For $x<y$ in $F$, the closure of $A_x\sm A_y$ is disjoint from $F_X$ so the set itself must be finite.
In other words, $\cA=\{A_x: x\in X\}$ is an almost chain of subsets of $N$. It is not difficult to check that (i) and (ii) hold.
\end{proof}

The characterization of $\propR$ via liftings in Lemma \ref{cde:1} translates very well to almost chains of subsets of $\omega$ which gives the following.

\begin{proposition}\label{cde:5}
Let $K=\ult(\fA/\fin)$ where $\fA$ is generated by an almost chain $\cA=\{A_x:x\in X\}$
of subsets of $\omega$ for some $X\sub (0,1)$.
The countable discrete extension $L=\ult(\fA)$ of $K$ has property $\propR$ if and only if
there is a family $\widetilde{\cA}=\{\widetilde{A_x} \sub \omega: x\in X\}$ such that
\begin{enumerate}[(i)]
\item $\widetilde{A_x}=^\ast A_x$ for every $x\in X$;
\item $\widetilde{A_x}\sub \widetilde{A_y}$ whenever $x,y\in X$ and $x<y$.
\end{enumerate}
\end{proposition}

\begin{proof}
We use Lemma \ref{cde:1}: To check that conditions are sufficient note that we can set
$\theta(A_x/\fin)=\widetilde{A_x}$
and extend $\theta$ to a lifting $\fA/\fin\to\fA$ since every $b\in \fA/\fin$ can be expressed as
a finite union of elements of the form $A_y/\fin - A_x/\fin$.
Necessity follows from the fact that, given a lifting $\theta: \fA/\fin\to\fA$,
the sets $\widetilde{A_x}=\theta(A_x/\fin)$ are as required.
\end{proof}

\section{Between $\propR$ and $\propE$} \label{between}

In this section we present a construction of a countable discrete extension of a compact line of weight $\omega_1$ without property $\propR$, but with an extension operator of norm $3$.
At the end of the section we will also apply this result to spaces which are not necessarily zero-dimensional.

The construction below and Theorem \ref{btw:2} is due to Witold Marciszewski. It will be convenient  to consider a subset $X$ of the Cantor set $2^\omega$ rather than of $[0, 1]$ and replace $\omega$ by $2^{<\omega}$.
We can do this as the space $2^\omega$ can be seen as a subset of the interval $[0, 1]$.

\begin{construction}[Marciszewski] \label{btw:1}
Let us consider the full dyadic tree $T=2^{<\omega}$. By $\preccurlyeq$ we denote the lexicographic order on $2^\omega\cup 2^{<\omega}$:
$x \preccurlyeq y$ means that either  $x$ is an initial segment of $y$ or $x(k) < y(k)$ for $k = \min \{n \in \omega : x(n) \neq y(n)\}$.

Take any set $X \subseteq 2^\omega$ such that for every
$x\in X$ both the sets $\{n: x(n)=0\}$ and $\{n: x(n)=1\}$ are infinite.
For each $x\in X$ we denote
\[ S_x = \{x|n ^\frown 0 : n \in \omega, x(n) = 1\}\]

Consider the family $\cA_X=\{A_x:x\in X\}$ where
\[A_x = \{t \in T: t \preccurlyeq x\} \sm S_x.\]
Note that $\cA_X$, a family of subsets of a countable set $T$, is an almost chain, that is $A_x\sub^\ast A_y$ whenever $x,y\in X$ and $x\preccurlyeq y$. Write $\fA_X$ for the Boolean algebra by $\cA_X \cup \fin(T)$.
\end{construction}

\begin{theorem}[Marciszewski]\label{btw:2}
In the setting of Construction \ref{btw:1}, if $X$ is uncountable, then the space $\ult(\fA_X/\fin(T))$ is a separable compact line
and the space $\ult(\fA_X)$ is its countable discrete extension
without property $\propR$.
\end{theorem}

\begin{proof}
In view of Lemma \ref{cde:3} we only have to show the lack of property $\propR$.

Suppose otherwise that $L=\ult(\fA_X)$ has property $\propR$; then by Lemma \ref{cde:5}, for every $x\in X$ there is a finite modification $C_x$ of $A_x$ so that the sets $C_x$ form a chain.
Consider  a function $\vf : X\to \omega$ defined by
\[\vf(x) = \min\{n \in \omega : C_x \triangle A_x \subseteq \bigcup_{j < n} 2^j\}.\]

Since the set $X$ is uncountable, there is some $k\in \omega$ such that
$Y = \vf^{-1}(\{k\})$ is also uncountable.
It follows that $Y$ has a left-sided accumulation point $y \in Y$, so there is a sequence $x_n \prec y$ in $Y$ such that $x_n\to y$.

As $y(m)=1$ infinitely often, there is some $m > k$ and $x\in Y$ satisfying
\[x|m = y|m,\ x(m) = 0 \text{ and } y(m) = 1.\]
Put $\sigma = x|(m+1)$. As $y(m) = 1$ and $\sigma(m) = 0$, we have $\sigma \in S_y$, which implies that $\sigma \notin A_y$. We also have $m>k$, so $\sigma \notin C_y$ (as $(C_y \triangle A_y) \cap 2^{m+1} = \emptyset$).

Now from the definition of $\sigma$ we have $\sigma \prec x$ and $\sigma \notin S_x$. This means that $\sigma \in A_x$. As $m > k$, it follows that $\sigma\in C_x$.
Finally, $\sigma \in C_x\sm C_y$, which means that sets $C_x$ do not form a chain, a contradiction.
\end{proof}

We can, however, construct a sequence of measures as in Lemma \ref{cde:2} to prove the following.

\begin{theorem}\label{btw:3}
In the setting of Construction \ref{btw:1},
if the set $X$ is uncountable, then the space $L = \ult(\fA_X)$ is a countable discrete extension
of $K = \ult(\fA_X/\fin(T))$ satisfying $\eta(K, L) = 3$.
\end{theorem}

\begin{proof}
For any $\sigma\in 2^{<\omega}$ we denote by $p(\sigma) \in \ult(\fA_X)$ the unique non-principal ultrafilter satisfying
\[A_x \in p(\sigma) \mbox{ if and only if } \sigma \prec x \mbox{  for } x\in X.\]

Given $n$ and  $\sigma\in 2^n$, if $\sigma(n-1) = 0$, then we denote by $\sigma'$ the sequence $\sigma|(n-1) ^\frown 1$.
Additionally, if we can find the greatest number $m < n-1$ such that $\sigma(m) = 0$, then we put $\sigma'' = \sigma|m ^\frown 1$.
Now we define
\begin{equation*}
\mu_\sigma =
\begin{cases}
\delta_{p(\sigma)} &\text{if } \sigma(n-1) = 1,\\
\delta_{p(\sigma)} - \delta_{p(\sigma')} + \delta_{p(\sigma'')} &\text{if } \sigma(n-1) = 0 \mbox{ and } \sigma|(n-1) \neq 1^{n-1},\\
\delta_{p(\sigma)} - \delta_{p(\sigma')} &\text{if } \sigma(n-1) = 0 \mbox{ and } \sigma|(n-1) = 1^{n-1}.
\end{cases}
\end{equation*}

\begin{scclaim}
The following hold for every $\sigma\in T$ and $x\in X$
\begin{enumerate}[(a)]
    \item if $\sigma\in A_x$ then  $ \mu_\sigma({A_x}) = 1$,
    \item if $\sigma\notin A_x$  then $ \mu_\sigma({A_x}) = 0$.
\end{enumerate}
\end{scclaim}

To verify part (a) of Claim, take $\sigma\in A_x\cap 2^n$. We have $\sigma \prec x$ and $\sigma\notin S_x$.
We know that $\delta_{p(\sigma)}({A_x}) = 1$, as $A_x\in p(\sigma)$ from the definition of $p(\sigma)$.
Now there are two cases.

\begin{itemize}
    \item[Case 1:] $\sigma(n-1) = 1$; then $\mu_\sigma = \delta_{p(\sigma)}$ and $\mu_\sigma({A_x}) = 1$.
    \item[Case 2:] $\sigma(n-1) = 0$; then either
    $\sigma|(n-1) \neq x|(n-1)$ which implies $ \sigma', \sigma'' \prec x$
    or $\sigma|(n-1) = x|(n-1)$ and  $\sigma', \sigma'' \succ x.$
    In either case,  $\delta_{p(\sigma')}(A_x) = \delta_{p(\sigma'')}(A_x)$ and
    \[\mu_\sigma({A_x}) = (\delta_{p(\sigma)} - \delta_{p(\sigma')} + \delta_{p(\sigma'')})({A_x}) = \delta_{p(\sigma)}(A_x) = 1.\]
\end{itemize}

To verify \textit{(b)} consider $\sigma\in 2^n$ such that $\sigma\notin A_x$; then $x \prec \sigma$ or $\sigma\in S_x$.

In the first case, we have $\delta_{p(\sigma)}({A_x}) = 0$ and similarly $\delta_{p(\sigma')}({A_x}) = 0, \delta_{p(\sigma'')}({A_x}) = 0$ (as $\sigma \preccurlyeq \sigma', \sigma''$ if they are defined).
It follows that $\mu_\sigma({A_x}) = 0$.

If $\sigma\in S_x$ then we have $\sigma \preccurlyeq x$ and $\sigma' \preccurlyeq x$, as $\sigma' = x|n$ and $\sigma \prec \sigma'$.
Moreover, if $\sigma''$ is defined then it is $\preccurlyeq$-bigger than $x$, as $\sigma|(n-1) = x|(n-1)$ and $\sigma|(n-1) \prec \sigma''$.
It follows that
\[\mu_\sigma({A_x}) = (\delta_{p(\sigma)} - \delta_{p(\sigma')} + \delta_{p(\sigma'')})({A_x}) = 1 - 1 + 0 = 0,\]
so the proof of Claim is complete.

The sequence $(\mu_\sigma)_{\sigma\in T}$ satisfies $\|\mu_\sigma\|\le 3$ and
Claim says that $\mu_\sigma-\delta_\sigma$ is zero on the elements of the generating chain of $\fA$.
It follows easily that the set
\[\{\sigma\in T: (\mu_\sigma-\delta_\sigma)(A)\neq 0\}\]
is finite for every $A\in\fA$, so our measures satisfy conditions of
Lemma \ref{cde:2} which means that $\eta(K,L)\le 3$.
On the other hand, $\eta(K,L)\ge 3$ by Theorem \ref{ca:1} and Theorem \ref{btw:2}.
\end{proof}

To conclude this section, we shall turn example from Theorem \ref{btw:3} into a general result.
Recall first the following fact concerning compact lines \cite[Theorem 2.4]{HL74}.

\begin{theorem} [Heath and Lutzer]\label{btw:4}
If  $K' \sub K$ are  compact lines then there is a norm-one extension operator $C(K')\to C(K)$.
\end{theorem}

Then recall the following observations --- the first one is (a slightly modified version of) \cite[Lemma 8.6]{MP18}.

\begin{lemma} \label{btw:5}
Let $K$ be a separable compact line of uncountable weight $\kappa$.
Then $K$ contains a topological copy of the space $2^\omega\!_X$, where $X$ is a dense subset of $2^\omega$ with $|X|=\kappa$.
\end{lemma}

\begin{remark}\label{btw:6}
If $K'\sub K$ are any compact spaces, then any countable extension $L'$ of $K'$ defines $L\in\cde(K)$ in an obvious way:
Say that $L'=K'\cup \omega$; then the space $L=K\cup\omega$ is obtained by considering a topological disjoint union of $K$ and $K'\cup \omega$ and identifying every point in $K'$ with its copy in $K$.
\end{remark}

\begin{theorem} \label{btw:7}
If $K$ is a nonmetrizable separable compact line, then there is $L \in \cde(K)$ satisfying $\eta(K, L) = 3$.
\end{theorem}

\begin{proof} As $K$ is nonmetrizable,  $\kappa=w(K) \geq \omega_1$.
By Lemma \ref{btw:5}, we can find inside $K$ a copy $K'$ of a zero-dimensional space
$2^\omega\!_X$, where $|X|=\kappa$.
Then, combining Theorem \ref{btw:3} with Lemma \ref{cde:4}, we can define $L'\in\cde(K')$ such that $\eta(K',L')=3$.
In turn, we get an `obvious' countable discrete extension $L=K\cup\omega$ mentioned in Remark \ref{btw:6}.

We know that there is an extension operator $E': C(K') \to C(L')$ of norm $3$, so we can define an extension operator $E: C(K) \to C(L)$ by
\begin{equation*}
Ef(x) =
\begin{cases}
f(x)          &\text{for } x \in K, \\
E'f|K' (x)    &\text{for } x \in L \sm K = \omega.
\end{cases}
\end{equation*}

Observe that $Ef$ is indeed continuous on $L$; clearly, $\|E\|=\|E'\|=3$.
On the other hand, by Theorem \ref{btw:4}, $\eta(K,L)<3$ would mean that $\eta(K', L') < 3$, which cannot hold.
\end{proof}

\section{Outside $\propE$}

Recall that if $\cI$ is a proper $\sigma$-ideal of subsets of, for instance,  the Cantor set $2^\omega$ then
\[  {\rm non}(\cI)=\min\{ |X|: X\notin \cI\}.\]
As in \cite{AMP20} we consider here the $\sigma$-ideal $\idE$ of subsets of $2^\omega$ that can be covered by a countable number of closed sets of measure zero.
Cardinal coefficients of $\idE$ are discussed by Bartoszy\'nski and Shelah \cite{BS92}; that ideal  is usually denoted by $\mathcal E$ --- we have changed the notation for an obvious reason.
Clearly, $\idE\sub \cN\cap\cM$
where $\cN$ denotes the family of $\lambda$-null subsets of $2^\omega$, where $\lambda$ is the standard product measure and $\cM$ denotes the $\sigma$-ideal of meager sets in $2^\omega$.
Hence
\[ \nonE \le \min\left( {\rm non}(\cN), {\rm non}(\cM)\right);\]
the strong inequality in the above formula is relatively consistent, see  \cite{BS92}.
Recall that cardinal coefficients of the classical $\sigma$-ideals do not change if we replace $2^\omega$ by any uncountable Polish space (and $\lambda$ by any nonatomic Borel measure on it), cf.\ \cite{BS92} and \cite{Fr5}.

The main point of this section is to prove the following theorem.

\begin{theorem}\label{out:2}
Let $\kappa \geq \nonE$. There is a zero-dimensional compact line $K$ of weight $\kappa$ which has a countable discrete extension without property $\propE$.
\end{theorem}

Let us first give a construction leading to the space mentioned in the theorem and recall some notions used in the proof.

\begin{construction} \label{out:1}
We consider a subtree $T$ of $\omega^{<\omega}$ defined as
\[T=\{\sigma: \sigma(n)\le n \mbox{ for every } n\},\]
and its body
\[ C=\{0\}\times \{0,1\}\times\{0,1,2\}\times\ldots\]
Again, as in \ref{btw:1}, we can consider the lexicographic order $\preccurlyeq$ on the set $T \cup C$. As the space $\omega^\omega$ can be identified with the set of irrational numbers, we can see $C$ as a subset of the unit interval $[0, 1]$.

Take any set $X \sub C$. For $x\in X$ we set
\[ S_x = \{x|n ^\frown 0 : x(n) \mbox{ is odd}\},\]
\[A_x = \{\sigma \in T: \sigma \preccurlyeq x\} \sm S_x\]
and consider the family $\cA_X=\{A_x:x\in X\}$.

It follows as in \ref{btw:1} that
$\cA_X$ is an almost chain of subsets of $T$.
Therefore, if we denote by $\fB_X$ the Boolean algebra generated by $\cA_X\cup\fin(T)$,
$K=\ult(\fB_X/\fin(T))$ is a separable compact line with $w(K) = |X|$ and
$L=\ult(\fB_X)$ is its countable discrete extension, see Lemma \ref{cde:3}.
\end{construction}

For the proof below, we equip $C$ with the standard product measure $\lambda$;
thus for every $\sigma\in T$ of length $n$ and $i\le n$ our measure satisfies
\begin{equation}\label{jeden}
\lambda([\sigma^\frown i])=\lambda([\sigma])/(n+1),
\end{equation}
where $[\sigma]$ denotes all elements of $C$ extending $\sigma$.
Then we may think that $\idE$ is the $\sigma$-ideal of subsets of the space $C$ generated by closed subsets of $\lambda$-measure zero.

We shall make use of the fact that the measure $\lambda$ satisfies the Lebesgue density theorem, that is for every closed set $F\sub C$ we have
\[\lim_{k\lra\infty} \frac{\lambda( [x|k]\cap F)]}{\lambda([x|k])}=1,\]
for $\lambda$-almost all $x\in F$. See e.g.\ \cite[Proposition 2.10]{Mi08} for a short proof that such a property is shared by every probability measure on a Polish ultrametric space.

\begin{proof} (of Theorem 5.2)
As $\kappa \geq \nonE$ let us fix any set $X\sub C$ of cardinality $\kappa$ and such that  $X\notin\idE$. This means that whenever $X=\bigcup_n X_n$ then $\lambda(\overline{X_n})>0$ for some $n$.

Take the space $L = \ult(\fB_X)$ from Construction \ref{out:1}. We shall check that $L$ is a countable discrete extension of $K = \ult(\fB_X/\fin(T))$ which does not have property $\propE$.
By Lemma \ref{cde:2} it is enough to demonstrate that whenever $(\mu_\sigma)_{\sigma \in T}$ is a sequence of finitely additive measures on $\fB_X$
satisfying
\begin{enumerate}[(i)]
  \item $\mu_\sigma(I)=0$ for every $I\in\fin(T)$ and every $\sigma$;
  \item the set $\{\sigma \in T: | \mu_\sigma(A) - \delta_\sigma(A)| \ge\eps\} $ is finite for every $A\in\fB_X$ and $\eps>0$,
\end{enumerate}
then $\sup_\sigma \|\mu_\sigma\|=\infty.$

Given such a family of measures $\mu_\sigma$, for every $x\in X$ we put
\[ B_x=\{\sigma\in T: |\mu_\sigma(A_x)-\delta_\sigma(A_x)|<1/4\}.\]
It follows that every $B_x$ is a finite modification of $A_x$ and we can write $X=\bigcup_n X_n$, where
$X_n$ is the set of those $x\in X$ for which $\sigma \in B_x$ is equivalent to $\sigma\in A_x$ for all $\sigma$ of length $\ge n$.

By the preparatory remarks, there is $n_0$ such that, writing $F=\overline{X_{n_0}}$ we have
$\lambda(F)>0$; in turn, there is a point $y\in F$ at which the set $F$ has density one.

Let us a fix a natural number $p$. Take $n>n_0$ such that
\begin{equation}\label{dwa}
 \frac{\lambda( [y|n]\cap F)]}{\lambda([y|n])}> \frac{2p}{2p+1}.
 \end{equation}
We can, of course, assume that $n+1=2pk$ for some natural number $k$. Consider the set
\[I=\{i\le n: y|n ^\frown i=x|(n+1)\mbox{ for some } x\in F\}.\]

\begin{scclaim}[Claim 1.]
The set $I$ satisfies\[ |I|\ge \frac{2p}{2p+1}(n+1).\]

Indeed, $i\notin I$ implies $[y|n ^\frown i]\cap F=\emptyset$, so the claim follows from (\ref{jeden}) and (\ref{dwa}).
\end{scclaim}

Divide $\{0,\ldots, n\}$ into $p$ consecutive intervals $J_0,\ldots, J_{p-1}$, each satisfying $|J_i|=2k$.

\begin{scclaim}[Claim 2.]
$|I\cap J_i|>k$ for every $i< p$.
\end{scclaim}

Indeed, by Claim 1 we have
\[ \frac{2p}{2p+1}(n+1)\le |I|=|I\cap J_i|+|I\sm J_i|\le |I\cap J_i|+n+1-2k, \mbox{ so}\]
\[ |I\cap J_i|\ge \frac{2p}{2p+1}(n+1)-(n+1)+2k=2k\frac{p+1}{2p+1}>k.\]
\medskip

It follows from Claim 2 that for every $i\le p$, the set $I\cap J_i$ contains at least one odd and one even number.
Pick an even number $m_0\in I\cap J_0$, odd $m_1\in I\cap J_1$, even $m_2\in I\cap J_2$ and so on.

Put $\sigma=y|n ^\frown 0$. From the definition of $I$ and the fact that $F=\overline{X_{n_0}}$
we conclude that there
are $x_0,\ldots, x_{p-1}\in X_{n_0}$ such that
\[ \sigma \in A_{x_0}, \sigma \notin A_{x_1}, \sigma \in A_{x_2}\ldots\]
Consequently, as $A_x$'s agree with $B_x$'s at that level,
\[ \mu_\sigma(A_{x_0})>3/4, \mu_\sigma(A_{x_1})<1/4, \mu_\sigma(A_{x_2})>3/4,\ldots \]
We know that $x_0\preccurlyeq x_1\preccurlyeq\ldots$, so the sets $A_{x_i}$ form an almost chain.
As the measure $\mu_\sigma$ vanishes on finite sets, we have, for instance,
\[ \mu_\sigma(A_{x_1}\sm A_{x_0})=\mu_\sigma(A_{x_1})-\mu_\sigma(A_{x_0})<1/4-3/4=-1/2.\]
In the same manner we get  $|\mu_\sigma|(A_{x_i}\sm A_{x_{i-1}})\ge 1/2$ and conclude that
 $\|\mu_\sigma\|\ge (p-1)/2$.

As $p$ was an arbitrary natural number, we get $\sup_\sigma\|\mu_\sigma\|=\infty$, as required.
\end{proof}

\begin{corollary} \label{out:3}
 Every separable compact line $K$ of weight $\ge\nonE$ has a countable discrete
 extension without property $\propE$.
\end{corollary}

\begin{proof}
We can argue as in the proof of Theorem \ref{btw:7}:
Find a zero-dimensional subspace $K'$ of $K$ and $L'\in\cde(K')$ such that $\eta(K',L')=\infty$, combine $L'$ with $K$ to obtain $L \in CDE(K)$ and note that $\eta(K,L)$ is also infinite by Theorem \ref{btw:4}.
\end{proof}

Let us also mention that a modification of the proof of Theorem \ref{btw:7}, that is thinning out the tree used there, should give $L\in\cde(K)$ with $\eta(K,L)$ finite but arbitrarily large.

\section{Remarks and problems}\label{final}

Let us briefly recall here that if $K$ is an arbitrary compact space and $L$ is a countable discrete extension of $K$ without property $\propE$, then one can form a nontrivial short exact sequence
\[ 0 \to c_0 \to C(L) \to C(K) \to 0,\]
which means that $c_0$ embeds onto an uncomplemented subspace $Z$ of $C(L)$ and the quotient $C(L)/Z$ is isomorphic to $C(K)$;
see \cite[Lemma 2.4]{MP18}.
Here $c_0$ is the classical Banach space of sequences converging to 0; its special role in  short exact sequences of Banach spaces stems from the fact that, by Sobczyk's theorem, $c_0$ is complemented in every separable Banach superspace.
Moreover, it was shown in \cite{AMP20} that for every nontrivial short exact sequence
\[ 0 \to c_0 \to \ ?  \to C(K) \to 0,\]
the Banach space in question can be found inside a space of the form $C(L)$, where $L$ is a countable discrete extension of $M_1(K)$, the dual unit ball in $C(K)^\ast$.
Twisted sums of Banach spaces and, more generally, homological methods in Banach space theory are a subject of a recent extensive monograph \cite{CC23}.

In the light of results from \cite{MP18}, \cite{CT19} and \cite{AMP20} the following problem seems to be natural and apparently remains open.

\begin{problem}\label{problem}
Is it relatively consistent that, whenever $K$ is a separable compact space of weight $\omega_1$, $\eta(K,L)<\infty$ for every $L\in \cde(K)$?
\end{problem}

In particular, having in mind the results presented above, one can ask if Problem \ref{problem} can be settled for the class of separable compact lines.

Recall that it is not very difficult to demonstrate that every \textit{nonseparable} compact line $K$ with $w(K)=\omega_1$ admits a countable discrete extension $L$ such that $\eta(K,L) = \infty$.
Namely, one can construct $L=K\cup\omega$ such that $L$ is the closure of the set $\omega\sub K$. Then, if we supposed that $\eta(K,L)<\infty$ then it would follow from Lemma \ref{p:1}(b) that $K$ must support a strictly positive nonnegative measure $\mu$.
However, a compact line carrying such a measure is necessarily separable
(see \cite[section 8]{MP18} for the details).

\end{document}